\documentclass{article}

\usepackage{amsmath}

\usepackage{amssymb}

\usepackage{graphicx}

\usepackage{color}

\usepackage[all]{xy}

\usepackage{amsthm}

\theoremstyle{plain}

\newtheorem{thm}{Theorem}

\newtheorem{lemma}{Lemma}

\newtheorem*{thm1}{Theorem}

\newtheorem*{lem1}{Lemma}

\newtheorem*{prop1}{Proposition}

\newtheorem*{cor1}{Corollary}

\newtheorem*{question}{Question}

\theoremstyle{remark}

\newtheorem*{claim}{Claim}

\newtheorem*{remark}{Remark}

\theoremstyle{definition}

\newtheorem*{dfn}{Definition}

\newcommand{\co}{\colon\thinspace}

\newcommand{\inv}[1]{#1^{-1}}

\newcommand{\bound}{\partial}

\newcommand{\AND}{\qquad \mathrm{and} \qquad}

\newcommand{\Q}{\mathbb{Q}}

\newcommand{\C}{\mathbb{C}}

\newcommand{\tr}{\mathrm{Tr}}

\definecolor{gray1}{gray}{0.6}

\begin{document}

\title{Trace fields and commensurability of link complements}

\author{E. Chesebro\footnote{Mathematics Department, Rice University, {\it chesebro@rice.edu}} and J. DeBlois\footnote{Department of Mathematics, Statistics and Computer Science, University of Illinois at Chicago, {\it jdeblois@math.uic.edu}}}

\maketitle

\begin{abstract}  This paper investigates the strength of the trace field as a commensurability invariant of hyperbolic 3-manifolds.  We construct an infinite family of two-component hyperbolic link complements which are pairwise incommensurable and have the same trace field, and infinitely many 1-cusped finite volume hyperbolic 3-manifolds with the same property.  We also show that the two-component link complements above have integral traces, but each has a mutant with a nonintegral trace.
\end{abstract}

\section{Introduction}

Manifolds $M$ and $M'$ are \emph{commensurable} if there is a manifold $N$ which is a finite cover for both $M$ and $M'$.  The study of the commensurability relation among hyperbolic knot and link complements in $S^3$ was initiated by W. Thurston, who gave examples of commensurable and incommensurable knots and links in Chapter 6 of his notes \cite{Th}.  Thurston's incommensurable examples may be distinguished using their \textit{trace fields}.  For a hyperbolic manifold $M = \mathbb{H}^3/\Gamma$, the trace field of $M$ is defined to be the field obtained by adjoining to $\mathbb{Q}$ the traces of elements of $\Gamma$.  Here $\Gamma$ is a torsion--free discrete subgroup of $\mathrm{PSL}_2(\mathbb{C})$.  Any two hyperbolic link complements in $S^3$ which are commensurable share a trace field (see Theorem 4.2.1 of \cite{MaR}).  The main theorem of this paper shows that the converse fails dramatically.

\begin{thm}  For any $k \geq 2$, there exist infinitely many commensurability classes of $k$-component links in $S^3$ each of whose complement has trace field $\mathbb{Q}(i,\sqrt{2})$.  \label{links} \end{thm}

\begin{figure}[ht]

\begin{center}

\includegraphics[height=1.25in]{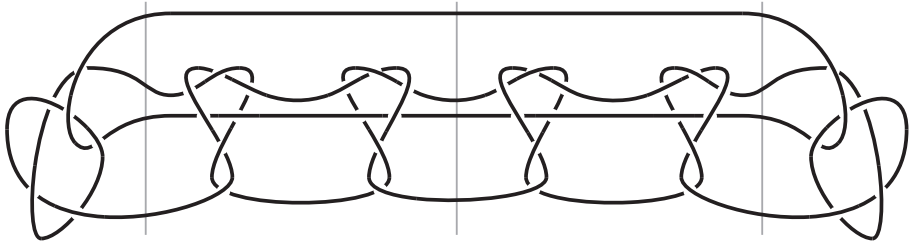}

\end{center}

\caption{The link $L_2$} \label{linkspic}

\end{figure}

Our primary object of study is a family of two-component links $L_n$, $n \in \mathbb{N}$, constructed from tangles $S \subset B^3$ and $T \subset S^2 \times I.$  The link $L_2$ is pictured in Figure \ref{linkspic}.  The gray vertical lines in the figure determine 4-punctured spheres which divide $L_2$, from left to right, into the tangle $S$, followed by two copies of $T$, capped off with the mirror image of $S$.  The link $L_n$ is constructed analogously, but with $n$ copies of $T$.

\bigskip

There is a unique hyperbolic structure on $B^3-S$ with totally geodesic boundary a 4-punctured sphere.  We describe a polyhedral decomposition for this in Section \ref{sec:S}.  In Section \ref{sec:T}, we give a polyhedral decomposition of a hyperbolic structure on $(S^2 \times I)-T$ with two totally geodesic 4-punctured sphere boundary components.  With these geometric structures, the boundary components of $(S^2 \times I)-T$ are isometric to the boundary of $B^3-S$, so that gluing copies of $B^3-S$ and $(S^2 \times I)-T$ together along their boundaries as in Figure \ref{linkspic} preserves the geometric structure of the pieces.  This allows us to give an explicit description, in Theorem \ref{Ln}, of the Kleinian group $\Gamma_n < \mathrm{PSL}_2(\mathbb{C})$ which uniformizes the hyperbolic link complement $M_n = S^3-L_n$.  This may be of independent interest.

Our description of the algebra and geometry of $M_n$ allows easy computation of many geometric invariants.  Among these, the \textit{complex modulus} of the cusps distinguishes the commensurability class of $M_m$ from that of $M_n$ for $m \neq n$.  The complex modulus of a cusp is a commensurability invariant of the Euclidean similarity class of a cross--section.  This is actually the invariant originally used by Thurston to distinguish the incommensurable examples in \cite{Th}; but as mentioned above, these examples also have different trace fields.  One might suppose that the complex modulus is simply a sharper commensurability invariant of cusped manifolds than the trace field, but even for hyperbolic knot complements this is not the case.  For example, the Figure 8 knot complement and the ``dodecahedral knot'' complements of Aitchison-Rubinstein \cite{AR} have commensurable cusp cross--sections, but the trace field in the first case is $\mathbb{Q}(\sqrt{-3})$ and in the second is $\mathbb{Q}(\sqrt{-3},\sqrt{5})$ (see \cite{NeR2}, Sections 9 and 10).

Another easy observation using our algebraic description of the $M_n$ is that the elements of the Kleinian groups $\Gamma_n$ all have integral traces.  On the other hand, we show in Section \ref{sec:nonint} that for each $n$, the complement of a certain mutant of $L_n$ has a nonintegral trace. Bass showed that if $M = \mathbb{H}^3/\Gamma$ where $\Gamma$ has an element with a nonintegral trace, there are closed essential surfaces in $M$ associated to this trace \cite{Ba}.  We say that such surfaces are \textit{detected by the trace ring}.  In our case, it is easily seen that for any $n \geq 2$, closed essential surfaces in the complements of both $L_n$ and its mutant can be obtained by attaching annuli parallel to the cusps to the punctured spheres separating copies of $T$.  It is interesting to note that although these surfaces are present in both families of links, the trace ring does not detect any closed surfaces in the $M_n$.  Integrality of traces is a commensurability invariant, and so for each $n$ we also obtain that the complement of $L_n$ is incommensurable with the complement of its mutant.

Our original motivation for this work was the following question, asked of us by Alan Reid.

\begin{question} Do there exist infinitely many hyperbolic knot complements in $S^3$ which share a trace field?

\end{question}

A negative answer to this question would imply a conjecture of Reid \cite{Reid}, that the commensurability class of any hyperbolic knot complement contains only finitely many others.  However, the results of this paper would seem to raise the probability of a positive answer.  Although we do not know if our techniques may furnish such an answer, we use them in Section \ref{sec:onecusp} to give one-cusped examples.

\begin{thm}  There exist infinitely many commensurability classes of one-cusped hyperbolic manifolds with invariant trace field $\mathbb{Q}(i,\sqrt{2})$. \label{onecusped} \end{thm}

Here we use the \textit{invariant} trace field, which in general is a subfield of the trace field, since for manifolds which are not link complements in $S^3$ the trace field may fail to be a commensurability invariant (cf. \cite{MaR}, Ch. 3).  The manifolds of this theorem are constructed using the same polyhedral building blocks as the links of Theorem \ref{links}, but with faces identified so that the totally geodesic boundaries of the resulting manifolds have full isometry group.  Their incommensurability again follows by considering the complex moduli of the cusps.

\bigskip

\noindent \textbf{\Large Acknowledgments}

\bigskip

The authors thank Ian Agol, Richard Kent, Chris Leininger, and Peter Shalen for helpful conversations, and Joe Masters for suggesting the complex modulus.  Thanks to Dick Canary for helping us with Lemma \ref{convexcore}, and special thanks to Alan Reid for suggesting these questions to us and for many helpful conversations and suggestions.


\section{The tangle $\mathbf{S}$}  \label{sec:S}

A description of $B^3-S$ as an identification space of an ideal octahedron seems to be well known, for instance it follows easily from results in \cite{ZP}, but we do not know of a reference for an explicit proof.  We give a detailed, self-contained proof here for completeness as well as to extract precise algebraic information.  First we describe an identification space $M_S$ of the ideal octahedron.  It is possible to convince oneself that $M_S$ is homeomorphic to $B^3-S$ by drawing pictures (analogous to Thurston's picture of the Tripus \cite{Th}), but our proof follows from a careful discussion of the relationship between the geometric, algebraic, and topological objects involved.

\begin{figure}[ht]

\setlength{\unitlength}{.1in}

\begin{picture}(40,25)

\put(3,0) {\includegraphics[height= 2.5in]{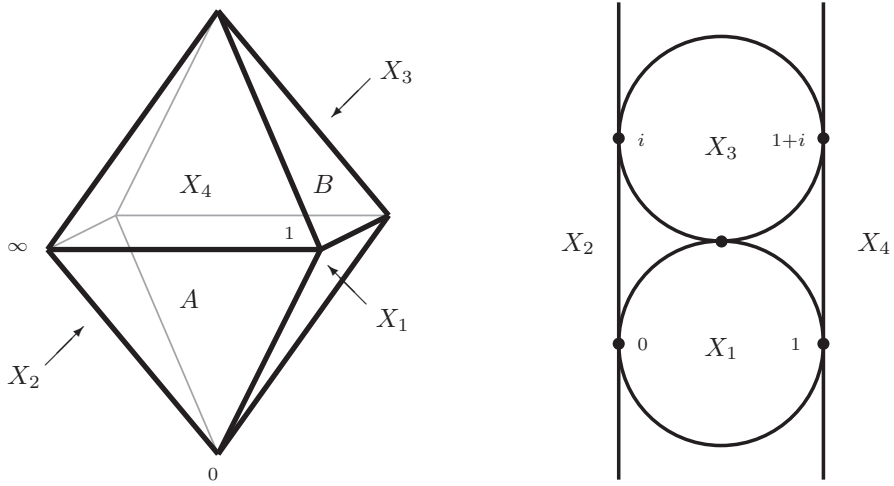}}

\put(17,15){$B$}

\put(10, 15){$X_4$}

\put(10, 9){$A$}

\put(19.8,9.2){\vector(-1,1){2}}

\put(20.3,8){$X_1$}

\put(1, 5){$X_2$}

\put(3,6){\vector(1,1){2}}

\put(20.5,21){$X_3$}

\put(20,21){\vector(-1,-1){2}}

\put(1,12){$\scriptstyle{\infty}$}

\put(15.5,12.6){$\scriptstyle{1}$}

\put(11.5,0){$\scriptstyle{0}$}

\put(30,12){$X_2$}

\put(45.5, 12){$X_4$}

\put(37.5, 6.5){$X_1$}

\put(37.5, 17){$X_3$}

\put(34,17.5){$\scriptstyle{i}$}

\put(34, 6.8){$\scriptstyle{0}$}

\put(41,17.5){$\scriptstyle{1+i}$}

\put(42, 6.8){$\scriptstyle{1}$}

\end{picture}

\caption{The ideal octahedron and a fundamental polyhedron $\mathcal{P}$ for $G$}  \label{cap}

\end{figure}

Let $M_S$ be the identification space of the hyperbolic regular ideal octahedron $\mathcal{O}$ (see Figure \ref{cap}) by face pairings $r$ and $s$, where $r$ takes $X_1$ to $X_2$ fixing the (ideal) vertex they share, and $s$ takes $X_3$ to $X_4$ so that the vertex they share goes to the vertex shared by $X_4$ and $X_2$.  Call the faces of $\mathcal{O}$ labeled $X_i$ the \textit{internal} faces and all others the \textit{external} faces.  Since each external face abuts internal faces on all sides with dihedral angles of $\pi/2$ and vice--versa, and the internal faces are identified in pairs, $M_S$ has totally geodesic boundary.

After positioning $\mathcal{O}$ in $\mathbb{H}^3$ so that the face $A$ shown has vertices at 0, 1, and $\infty$, and the other ideal points have positive imaginary part, $\mathcal{O}$ is a subpolyhedron of the polyhedron $\mathcal{P}$ indicated on the right in Figure \ref{cap}.   The lines and circles of the figure lie in $\bound \mathbb{H}^3 = \C \cup \{\infty \}$ and bound totally geodesic hyperplanes which contain the faces $X_1,$ $X_2,$ $X_3,$ and $X_4$ of $\mathcal{O}.$  Each line and circle is labeled by the face of $\mathcal{O}$ which is contained by the corresponding hyperplane.  Each hyperplane divides $\mathbb{H}^3$ into two half-spaces, one of which contains the remaining hyperplanes.  The polyhedron $\mathcal{P}$ is the intersection of all these half-spaces.  The heavy dots, together with $\infty$, are ideal vertices of $\mathcal{P}$.  $\mathcal{O}$ is recovered by truncating the open triangular ends of $\mathcal{P}$ by geodesic hyperplanes perpendicular to its faces.  For instance, $A$ is obtained by truncating the end bounded by $X_1$, $X_2$, and $X_4$ by the geodesic hyperplane $\mathcal{H}$ with boundary $\mathbb{R} \cup \{\infty\}$.

The face pairings $r$ and $s$ are now realized by the isometries

\[ r= \left(\begin{array}{cc} 1 & 0 \\ -1 & 1 \end{array}\right) \AND
   s= \left(\begin{array}{cc} 2i & 2-i \\ i & 1-i \end{array}\right) \]
which give corresponding face parings on the polyhedron $\mathcal{P}.$  Let $G= \langle r,s \rangle < \mathrm{PSL}_2(\mathbb{C})$.   According to the Poincar\'{e} polyhedron theorem (see for example \cite{Ra}, Theorem 11.2.2), the quotient of $\mathbb{H}^3$ by $G$ is isometric to the identification space of $\mathcal{P}$ by the side pairings induced by $r$ and $s$ and has $G$ as its fundamental group. Since $r$ and $s$ also induce the side pairings of the internal faces of $\mathcal{O}$ prescribed to produce $M_S$, the inclusion $\mathcal{O} \hookrightarrow \mathcal{P}$ induces an isometric embedding $i_S \co M_S \rightarrow \mathbb{H}^3/G$.

 Recall that, for a Kleinian group $\Gamma$, the {\it convex core} of $M=\mathbb{H}^3/\Gamma$, denoted $C(\Gamma)$, is the smallest convex submanifold of $M$ with the property that the inclusion induced homomorphism between fundamental groups is an isomorphism.

\begin{lemma}  The image of $i_S$ is $C(G)$. \label{M_SandG} \end{lemma}

\begin{proof}

The complement of $\mathcal{O}$ in $\mathcal{P}$ has four components, each of which is homeomorphic to the product of an external face with $(0,\infty)$, with second coordinate given by distance to the face.  The action of $G$ preserves this stratification and identifies the faces of these components corresponding to the edge identifications of the external faces of $\mathcal{O}$.  Accordingly, we have that $\mathbb{H}^3/G$ is naturally homeomorphic to 

\[ M_S \cup_{\partial M_S} \left(\partial M_S \times [0,\infty)\right) \]
Therefore $G \cong \pi_1 M_S$.

Let $\widetilde{M}_S$ be the universal cover of $M_S$.  The developing map $\widetilde{M}_S \rightarrow \mathbb{H}^3$ determined by our embedding of the octahedron has image

\[ \bigcup_{g \in G} g.\mathcal{O}. \]

This follows from the polyhedron theorem, which asserts that the universal cover of $\mathcal{P}/G$ has developing image $\mathbb{H}^3$, tiled by $G$--translates of $\mathcal{P}$.  It now follows from work of Kojima (\cite{Ko1}, see also \cite{Ko2}) that the convex hull of $G$ has boundary consisting of translates of the hyperplanes containing external faces of $\mathcal{O}$.  These cover $\partial M_S$, and so $i_S(M_S)$ is the convex core of $\mathbb{H}^3/G$.

\end{proof}

\begin{figure}[ht]

\begin{center}

\includegraphics[width=3in]{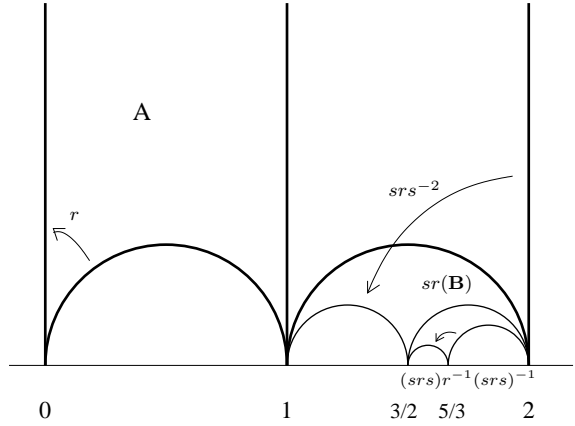}

\caption{The totally geodesic boundary of $M_S$}

\label{Sface}

\end{center}

\end{figure}

A combinatorial analysis of the edge identifications of external faces of $\mathcal{O}$ induced by $G$ reveals that $M_S$ has a single 4-punctured sphere boundary component $F$.  As mentioned above, $F$ is totally geodesic.  Let $\mathcal{H}$ be the hyperplane in $\mathbb{H}^3$ with $\partial \mathcal{H} = \mathbb{R} \cup \{\infty\}.$  Since $A \subset \mathcal{H}$ there is a subgroup $\Lambda < \mathrm{PSL}_2(\mathbb{R})$, stabilizing $\mathcal{H}$, which corresponds to $F.$  A fundamental domain for the action of $\Lambda$ is pictured in Figure \ref{Sface}, along with a set of face--pairing isometries generating $\Lambda$.  

An all--parabolic generating set for $\Lambda$ is given by

\[ \begin{array}{rcccl}

p_1 & = & r^{-1} & = &  \begin{pmatrix} 1 & 0 \\ 1 & 1 \end{pmatrix}  \\ 

p_2 & = &  rsrs^{-2} & = & \begin{pmatrix} -1 & 5 \\ 0 & -1 \end{pmatrix}  \\ 

p_3 & = & (srs)r^{-1}(srs)^{-1} & = & \begin{pmatrix} -14 & 25 \\ -9 & 16 \end{pmatrix}. \end{array} \]
The final conjugacy class of parabolic elements in $\Lambda$ is represented by 

$$p_4 \ = \  p_1 p_2 p_3^{-1} \ = \  \begin{pmatrix} 29 & -45 \\ 20 & -31 \end{pmatrix}.$$

\bigskip 

Evidently $p_1$ and $p_3$ are conjugate in $G$; so too are $p_2$ and $p_4$.  Combinatorial examination of the identification space reveals that $M_S$ has two rank one cusps, and so every parabolic element of $G$ is conjugate to one of $p_1$ or $p_2$.

The lemma below clarifies the relationship between $B^3-S$ and $G$.

\begin{figure}[ht]

\setlength{\unitlength}{.1in}

\begin{picture}(40,15)

\put(18,0) {\includegraphics[height= 1.45in]{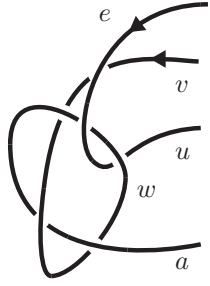}}

\put(22.75,13.5){$e$}

\put(24.75,4.25){$w$}

\put(26.75, .5){$a$}

\put(26.75, 9.75){$v$}

\put(26.75, 6.25){$u$}

\end{picture}

\caption{The tangle $S$ with labeled Wirtinger generators}  \label{S}

\end{figure}

\begin{lemma} With $S$ as in Figure \ref{S}, take the base point for $\pi_1(B^3-S)$ to be on the boundary surface high above the projection plane, and let Wirtinger generators correspond in the usual way to labeled arcs of the figure.   There is a faithful representation $\rho \co \pi_1(B^3-S) \rightarrow G$ with $\rho(a) = p_1$, $\rho(e) = p_2$, and $\rho(v) = p_3^{-1}$. 

\label{Srep} \end{lemma}

\begin{proof}

Reducing a standard Wirtinger presentation for $\pi_1(B^3-S)$, we obtain 

$$\left \langle a, w, e \,\Big| \, ew\inv{e}a = awa\inv{w} \right \rangle$$

Consider the following sequence of Tietze transformations.

$$\begin{array}{c}

\left \langle a, w, e, b \, \Big| \, b=\inv{e}aw, \ ew\inv{e}a = awa\inv{w} \right \rangle \\ \\

\left \langle a, w, e, b \, \Big| \, w=\inv{a}eb, \ ew\inv{e}a = awa\inv{w} \right \rangle \\ \\

\left \langle a, e, b \, \Big| \,  e\inv{a}eb\inv{e}a = a\inv{a}eba\inv{b}\inv{e}a \right \rangle \\ \\

\left \langle a, e, b \, \Big| \,  \inv{a}eb = ba\inv{b} \right \rangle \\ \\

\left \langle a, e, b \, \Big| \,  e = abab^{-2} \right \rangle \\ \\

\left \langle a, b \right \rangle

\end{array}$$

Recall that $G$ is generated by $r$ and $s$.  In fact, $G$ is {\it freely} generated by $r$ and $s$.  This follows from a standard ping--pong argument, since the sides of the fundamental polyhedron for $G$ do not intersect in $\mathbb{H}^3$.  Hence, the representation $\rho \co \pi_1(B^3-S) \longrightarrow G \subset \mathrm{PSL}_2(\mathbb{C})$ given by $a \mapsto r$ and $b \mapsto s$ is an isomorphism onto $G.$  Notice that the subgroup of $\pi_1(B^3-S)$ corresponding to the 4-punctured sphere $\bound B^3-\bound S$ is freely generated by $a,$ $v,$ and $e.$  It is easily checked that

$$\begin{array}{ccccc}

\rho(v) & = & \begin{pmatrix} 16 & -25 \\ 9 & -14 \end{pmatrix}  & = & p_3^{-1} \\ \\

&& \mathrm{and} \\ \\

\rho(e) & = & \begin{pmatrix} -1 & 5 \\ 0 & -1 \end{pmatrix} & = & p_2 .

\end{array}$$
Therefore $\rho$ takes $\pi_1(\partial B^3 - S)$ isomorphically to $\Lambda$.  The fact that meridians are taken to parabolic elements follows immediately from the observation that any meridian of $S$ is conjugate in $\pi_1(B^3-S)$ to either $a$ or $e$, and these are taken to $p_1$ and $p_2$ respectively.

\end{proof}

We wish to show that $M_S$ gives $B^3-S$ a hyperbolic structure with totally geodesic boundary.  The terminology we use to make this precise is in J. Morgan's paper on Thurston's proof of geometrization for Haken manifolds \cite{Mo}, for example.  Let $N(S)$ be a small open tubular neighborhood of $S$ in $B^3$.  Then $B^3-N(S)$ is a compact manifold with genus two boundary, and two annuli on the boundary defined by $\partial N(S)$.  The pair $(B^3 - N(S),\partial N(S))$ is a {\it pared manifold} (\cite{Mo}, page 58), and $(B^3-N(S))-\partial N(S)$ is homeomorphic to $B^3-S$.  When $\Gamma$ is a Kleinian group we write $C(\Gamma)$ to denote the convex core of $\mathbb{H}^3/\Gamma$.

\begin{lemma}  Let $(M,P)$ be a pared manifold, and suppose that $\rho: \pi_1 M \rightarrow \Gamma < \mathrm{PSL}_2 (\mathbb{C})$ is a faithful representation onto a non-Fuchsian geometrically finite Kleinian group where $C(\Gamma)$ has totally geodesic boundary.  If there is a one--to--one correspondence under $\rho$ between conjugacy classes of $\pi_1(M)$ corresponding to elements of $P$ and conjugacy classes of maximal parabolic subgroups of $\Gamma$, then $\rho$ is induced by a homeomorphism of $M - P$ to $C(\Gamma)$.  \label{convexcore} \end{lemma}

This is a ``homotopy hyperbolic manifolds are hyperbolic'' lemma for pared manifolds.  While it seems well known to experts in Kleinian groups, we do not know of a reference for a written proof, and it seems worth writing down as it may fail dramatically if $C(\Gamma)$ does not have totally geodesic boundary (see \cite{CM} for a thorough exploration of this phenomenon).  The proof is not difficult, following easily from results in \cite{CM} for example, but it requires introduction of the characteristic submanifold machinery.  Since this is somewhat outside the scope of the rest of the paper, we defer the proof to an appendix.

With $M = B^3 - N(S)$, $P = \partial N(S)$, $\Gamma = G$, and $\rho$ as in Lemma \ref{Srep}, it is clear that the hypotheses of Lemma \ref{convexcore} are satisfied .  Thus $B^3 - S$ is homeomorphic to $C(G)$, which by Lemma \ref{M_SandG} is isometric to the identification space $M_S$.  Call the induced homeomorphism $h_S \co B^3 - S \rightarrow M_S$.


\section{The tangle $\mathbf{T}$} \label{sec:T}

In this section we give a description of $(S^2 \times I) - T$ as an identification space of ideal polyhedra yielding a manifold with two totally geodesic boundary components.  In fact we show that each totally geodesic boundary component is isometric to the totally geodesic boundary of $B^3 - S$.  As far as we know, the description of the hyperbolic structure on $(S^2 \times I) - T$ was not previously known.  The method of proof is analogous to that of the previous section, and in some cases we skip details that were covered before.

Observe that $S^2 \times I$ has a mirror symmetry preserving $T$ and exchanging its boundary components.  Let $T_0$ be the tangle in $S^2 \times I$ pictured in Figure \ref{T0}, which is obtained from $T$ by cutting along the mirror locus of this reflection.

\begin{figure}[ht]

\setlength{\unitlength}{.1in}

\begin{picture}(40,25)

\put(3.5,0) {\includegraphics[height= 2.5in]{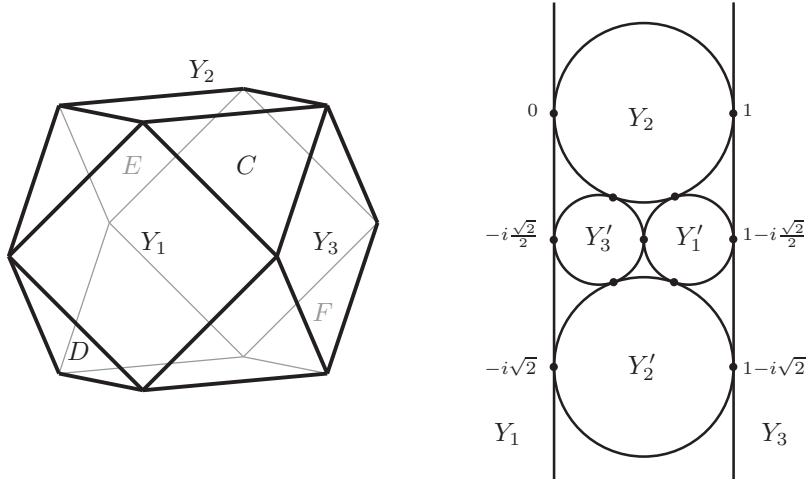}}

\put(15.5,16){$C$}

\put(10.5, 12){$Y_1$}

\put(6.65, 6.25){$D$}

\put(9.5,16){$\textcolor{gray1}{E}$}

\put(19.5,8.25){$\textcolor{gray1}{F}$}

\put(13,21){$Y_2$}

\put(19.5,12){$Y_3$}

\put(43, 2){$Y_3$}

\put(29, 2){$Y_1$}

\put(36, 18.5){$Y_2$}

\put(36, 5.5){$Y'_2$}

\put(33.75,12.25){$Y'_3$}

\put(38.5,12.25){$Y'_1$}

\put(28.5,12.5){$\scriptstyle{-i\frac{\sqrt{2}}{2}}$}

\put(42,19){$\scriptstyle{1}$}

\put(28.5, 5.5){$\scriptstyle{-i\sqrt{2}}$}

\put(30.75,19){$\scriptstyle{0}$}

\put(42, 5.5){$\scriptstyle{1-i\sqrt{2}}$}

\put(42,12.5){$\scriptstyle{1-i\frac{\sqrt{2}}{2}}$}

\end{picture}

\caption{The ideal cuboctahedron and a fundamental polyhedron for $H_0$}  \label{middle}

\end{figure}

\begin{figure}[ht]

\setlength{\unitlength}{.1in}

\begin{picture}(40,15)

\put(18,0) {\includegraphics[height= 1.45in]{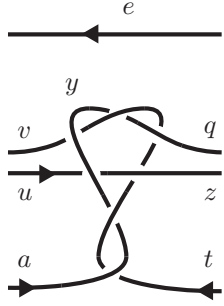}}

\put(24,15){$e$}

\put(21,11){$y$}

\put(28.25, 1.75){$t$}

\put(28.25, 8.75){$q$}

\put(28.25, 5.25){$z$}

\put(18.5, 1.75){$a$}

\put(18.5, 8.5){$v$}

\put(18.5, 5.25){$u$}

\end{picture}

\caption{The tangle $T_0$ with labeled Wirtinger generators}  \label{T0}

\end{figure}

We describe an identification space $M_{T_0}$ of the regular ideal cuboctahedron, the polyhedron pictured on the left in Figure \ref{middle}.  In the figure, interpret the square face opposite that labeled $Y_i$ to be labeled $Y_i'$.  Let $f$ be the parabolic taking $Y_2$ to $Y_1$, fixing the ideal vertex they share.  Let $g$ be the parabolic taking $Y_3$ to $Y'_1$, fixing the ideal vertex they share.  Let $h$ be the composition of the elliptic rotating $Y_3'$ by $\pi$ about its center with the parabolic taking $Y'_2$ to $Y'_3$ fixing the vertex they share.

Position the ideal cuboctahedron so that the face $C$ has ideal points at 0, 1, and $\infty$, where the vertex shared by $Y_1$ and $Y_2$ is at 0 and the vertex shared by $Y_1$ and $Y_3$ is at infinity, and all other ideal points of the cuboctahedron have negative imaginary part.  Then the face pairing isometries described above are realized by

\begin{align*} 
  & f\ =\ \begin{pmatrix} 1 & 0 \\ -1 & 1 \end{pmatrix} &
  & g\ =\ \begin{pmatrix} -1+i\sqrt{2} &  1-2i\sqrt{2} \\ -2 & 3-i\sqrt{2} \end{pmatrix}  \end{align*}

  \begin{align*} & h\ =\ \begin{pmatrix} 2i\sqrt{2} & -3-i\sqrt{2} \\ -3+i\sqrt{2} & -3i\sqrt{2} \end{pmatrix}. 
\end{align*}

Let $H_0= \langle\, f,g,h\, \rangle < \mathrm{PSL}_2(\mathbb{C})$.   A fundamental polyhedron for $H_0$ is indicated on the right in Figure \ref{middle}.  As with the ideal octahedron, the cuboctahedron is recovered from the fundamental polyhedron for $H$ by truncating open ends by perpendicular hyperplanes.  The eight resulting triangular faces are external, and the six square faces are internal.  Again, the regular ideal cuboctahedron has all right dihedral angles, hence the boundary of $M_{T_0}$ is totally geodesic.  We have

\begin{lemma}  The inclusion of the cuboctahedron into the fundamental polyhedron for $H_0$ induces an isometry $i_{T_0}$ from $M_{T_0}$ to $C(H_0)$.  \label{M_T0andH0}  \end{lemma}

\begin{proof} The proof is analogous to the proof of Lemma \ref{M_SandG}.  \end{proof}

A combinatorial analysis of the identification space $M_{T_0}$ reveals that it has two totally geodesic boundary components, each composed of four triangular faces and homeomorphic to a 4--punctured sphere.  Let $F_0$ and $F_1$ be these two totally geodesic punctured spheres, labeled so that $F_0$ contains the face $C$ of the cuboctahedron.  $F_0$ and the boundary $F$ of $M_S$ are combinatorially identical; this is reflected by the fact that $\Lambda$ is the boundary subgroup of $H_0$ fixing the totally geodesic hyperplane $\mathcal{H}$.  Indeed, we have

\[ \begin{array}{rcccc}

p_1 & = & r^{-1} & = & f^{-1} \\ 

p_2 & = &  rsrs^{-2} & = & fg^{-1}f^{-1}h^{-1}g \\ 

p_3 & = & (srs)r^{-1}(srs)^{-1} & = & (g^{-1}f^{-1}h)g^{-1}(h^{-1}fg). 

\end{array} \]
This shows that the convex cores $C(G)$ and $C(H_0)$ have isometric boundary components, each covered by $\mathcal{H}$.  Using the isometric embeddings $i_S$ and $i_{T_0}$ gives an isometry $i_0 \co F \rightarrow F_0.$

Observe that $H_0$ is normalized by 

$$\sigma \ = \ \begin{pmatrix} i & i-\sqrt{2} \\ 0 & -i \end{pmatrix}$$
whose action on the generators $f,$ $g,$ and $h$ is given by \begin{align*}
& \sigma f \sigma^{-1} = fgf^{-1} \\
& \sigma g \sigma^{-1} = (fg^{-1})f(gf^{-1}) \\
& \sigma h \sigma^{-1} = (fg^{-1})h^{-1}(gf^{-1}). \end{align*}

Conjugation by $\sigma$ induces an order-2 orientation preserving isometry on $M_{T_0}$.  The elements $\sigma p_i \sigma^{-1},$ $i \in \{1,2,3\},$ all preserve the geodesic hyperplane in $\mathbb{H}^3$ over the line $\mathbb{R} - i\sqrt{2}$.  This hyperplane contains a triangular face of the cuboctahedron projecting to $F_1 \subset M_{T_0}$.  Therefore the isometry induced by $\sigma$ interchanges the boundary components $F_0$ and $F_1$ of $M_{T_0}$.

\begin{lemma}  With $T_0$ as in Figure \ref{T0}, take the base point for $\pi_1 (S^2 \times I) - T_0$ to lie on the boundary component corresponding to endpoints $a$, $u$, $v$ and $e$, high above the projection plane, and take Wirtinger generators in the usual way.  There is a faithful representation $\pi_1 (S^2 \times I) - T_0 \rightarrow H_0$, which maps $a$, $e$, and $v$ to $p_1$, $p_2$, and $p_3^{-1}$, respectively.

\label{T0rep} \end{lemma}

\begin{proof}  $(S^2 \times I)-N(T_0)$ may be isotoped in $S^3$ to a standard embedding of a genus-3 handlebody.  Thus $\pi_1((S^2 \times I)-T_0)$ is free on three generators.  We claim that the group is generated by $a$, $e$, and $t$.  This follows after noticing that $v=\inv{y}xy$ where $y=\inv{(ta)}a(ta)$ and $x=\inv{(azq)}t(azq)=\inv{(ate)}t(ate)$.   (The relation $zq=te$ used in the last equality comes from the relationship between four peripheral elements in a 4-punctured sphere group.)  So far, we have established that $v, y \in \langle a, e, t \rangle.$  Now using the other punctured sphere relation, we have $u = a^{-1}ev \in \langle a, e, t \rangle.$  Finally, $z=yuy^{-1}$ and $q=z^{-1}te.$  Therefore $a,$ $e,$ and $t$ generate the group as claimed.

As in the previous section, it is easily seen that $H_0$ is freely generated by $f$, $g$, and $h$.  For our purposes, a more convenient free generating set for $H_0$ is 

$$\{f,\, fgf^{-1},\, p_2  \}.$$ 
Note that all of these generators are parabolic and peripheral, and $\sigma$ acts on them by interchanging the first two and taking the third to its inverse.  The representation of $\pi_1((S^2 \times I) - T_0)$ given by \begin{align*}
& a \mapsto f & & t \mapsto fgf^{-1} & & e \mapsto p_2
\end{align*}
is clearly faithful, and it is easily checked that $v$ maps to $p_3^{-1}$.  Because $u = a^{-1}ev$ is mapped to $p_1p_2p_3^{-1} = p_4$, we conclude that meridians are mapped to parabolic elements and that $\pi_1 ((S^2 \times \{0\})-T_0)$ is taken to $\Lambda$.

There is a visible involution of $(S^2 \times I)-T_0$ given by a rotation of $\pi$ around a vertical line down the center of the diagram in Figure \ref{T0}.  This involution exchanges the two boundary components.  With a proper choice of path between our basepoint and its image under this involution, the corresponding action on $\pi_1((S^2 \times I)-T_0)$ is given by \begin{align*}
  & a \leftrightarrow t & & e \leftrightarrow e^{-1} \end{align*}
This is the pullback of the aforementioned action by $\sigma$ on $H_0.$  Hence our representation maps $\pi_1 ((S^2 \times \{1\})-T_0)$ to $\sigma \Lambda \sigma^{-1}$ .

\end{proof}

As in the previous section, it follows from Lemma \ref{convexcore} and the lemmas above that $(S^2 \times I) - T_0$, the convex core $C(H_0),$ and $M_{T_0}$ are all homeomorphic.  Call the homeomorphism $h_{T_0} \co (S^2 \times I) - T_0 \rightarrow M_{T_0}$, and note that $h_{T_0}$ takes $(S^2 \times \{i\})-T_0$ to $F_i$ for $i = 0$ or 1.

We obtain a new hyperbolic manifold $M_{T}$ by doubling $M_{T_0}$ across the boundary component $F_1$.  Using the homeomorphism $h_{T_0}$, we see that the double of $(S^2 \times I) - T_0$ across $(S^2 \times \{1\}) - T_0$ is homeomorphic to the double of $M_{T_0}$ across $F_1$.  That is, we have a homeomorphism $h_T \co (S^2 \times I) - T \rightarrow M_T$.

In order to describe the Kleinian group whose convex core is isometric to $M_T$, we introduce some notation.  If $\Theta$ is a subgroup of a group $\Gamma$ and $g \in \Gamma$ then we write $\Theta^g$ for $g \Theta\inv{g}$.  If $N$ is an oriented manifold then $\overline{N}$ is the manifold obtained from $N$ by reversing the orientation.  If $\Gamma$ is a subgroup of $\mathrm{PSL}_2(\C)$ we denote the group whose elements are complex conjugates of the elements of $\Gamma$ by $\overline{\Gamma}$.  Similarly, if $g \in \mathrm{PSL}_2(\C)$ then $\overline{g}$ is the complex conjugate of $g$.  Note that complex conjugation induces the isometry of $\mathbb{H}^3$ given by reflection through $\mathcal{H}$, and that $\overline{\Gamma}$ is the conjugate of $\Gamma$ by this isometry.  Hence complex conjugation induces an orientation--reversing isometry between convex cores $C(\Gamma)$ and $C(\overline{\Gamma})$.

Now consider the element

 \[ c \ = \ \begin{pmatrix} 1 & i\sqrt{2} \\ 0 & 1 \end{pmatrix}  \ \in \ \mathrm{PSL}_2(\mathbb{C}) \]
and observe that $c$ takes the hyperplane over $\mathbb{R}-i\sqrt{2}$ to $\mathcal{H}$.  Recall that $F_1$ is the quotient of the hyperplane over $\mathbb{R}-i\sqrt{2}$ by $\Lambda^\sigma$.  Therefore in $H^c_0$ the corresponding surface subgroup $\Lambda^{c \sigma}$ fixes $\mathcal{H}$ and $\overline{H^c_0} \cap H^c_0 = \Lambda^{c\sigma}$.

We are in a position to apply the following lemma, which combines Maskit's combination theorem with observations of Morgan (compare with the statement of Theorem 8.2 of \cite{Mo}).

\begin{lemma}[Maskit, Morgan]  Suppose $G_0$ and $G_1$ are Kleinian groups with $\Lambda = G_0 \cap G_1$ a totally geodesic surface subgroup whose fixed circle divides $S^2$ into $B_0$ and $B_1$, with the property that for each $i$, the only elements $g_i \in G_i$ such that $g_i B_{i+1} \cap B_{i+1} \neq \emptyset$ are those in $\Lambda$.  Then $G = \langle\, G_0, G_1\, \rangle < \mathrm{PSL}_2(\mathbb{C})$ is a Kleinian group and as an abstract group we have
\[ G \cong G_0 *_{\Lambda} G_1. \]
The convex core $C(G)$ is isometric to $C(G_0) \cup_{F} C(G_1)$, where $F$ is the boundary component of $C(G_i)$ corresponding to $\Lambda$.  \label{Maskit}

\end{lemma}

\begin{proof}  The entire lemma except for the last sentence is the statement of Maskit's combination theorem for free products with amalgamation \cite{Mask}, in the special case of amalgamation along a totally geodesic surface.  The final sentence follows from Morgan's analysis below Theorem 8.2 of \cite{Mo}, noting that in this case the surface $X = f^{-1}(\frac{1}{2})$ \textit{is} the totally geodesic surface uniformized by $\Lambda$.

\end{proof}

Lemma \ref{Maskit} and the discussion above it show that $M_T$ is isometric to the convex core of $\langle\, H^c_0,\overline{H^c_0} \,\rangle$.  For our purposes it will be more convenient to use a conjugate of this group by $c$. Define

\[  H = \langle\, H^{c}_0, (\overline{H^c_0})\,\rangle^c = \langle\, H^{c^2}_0, \overline{H}_0\, \rangle. \]
(Note that $\bar{c} = c^{-1}$, which allows the simplified description of $H$ above.)  Clearly, the convex core $C(\langle\, H^c_0,\overline{H^c_0} \,\rangle)$ is isometric to the convex core $C(H)$.  This gives an isometry $i_T \co M_T \rightarrow C(H).$


\section{Main theorem} \label{sec:main}

Recall from the introduction that our main object of study is the link $L_n$, constructed by adjoining $n$ copies of $T$ end-to-end and capping off on the left by $S$ and the right by $\overline{S}$.  Let $M_n = S^3 - L_n$.  In this section we give an explicit algebraic and combinatorial description of $M_n$ and use this to prove our main theorem.

\begin{thm}  For each $n \in \mathbb{N}$, $M_n$ is hyperbolic and homeomorphic to $\mathbb{H}^3/\Gamma_n$, where

\[ \Gamma_n \ = \ \left\langle\, G^{c^{2n}},H^{c^{2(n-1)}}, \hdots , H, \bar{G}\, \right\rangle. \]

\label{Ln} \end{thm}

Before we prove the theorem, we set some notation and state and prove a lemma.  Let $j_0 \co \partial B^3 - S \rightarrow (S^2 \times \{0\}) - T_0$ be the gluing homeomorphism of the 4-punctured sphere which produces the complement of the tangle $S \cup T_0$ of Figure 7 from the tangles $S$ and $T_0$.

\bigskip

\begin{figure}[ht]

\begin{center}

\includegraphics[height=1.25in]{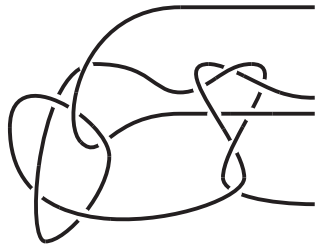}

\end{center}

\caption{$S \cup T_0$} \label{tanglepic}

\end{figure}

With our choice of base points, the map on fundamental groups induced by $j_0$ identifies the generators $a$, $u$, and $v$ as labeled in Figures \ref{S} and \ref{T0}.  Because the representations to $G$ and $H_0$ constructed in Lemmas \ref{Srep} and \ref{T0rep} identify these same generators in $\Lambda$, the diagram below commutes at the level of $\pi_1$ and so commutes up to isotopy. ($i_0$ is the isometry defined below Lemma \ref{M_T0andH0}.)
\[ \xymatrix{  F \ar[r]^{i_0} & F_0 \\  \partial B^3 - S \ar[u]^{h_S} \ar[r]^{j_0} & (S^2 \times \{0\}) - T_0 \ar[u]^{h_{T_0}} }  \]
Let $i_1 \co F_0 \rightarrow F_1$ be the restriction of the mirror symmetry of $M_T$ to $F_0$ and let $j_1 \co (S^2 \times \{0\}) - T \rightarrow (S^2 \times \{1\}) - T$ be the homeomorphism induced by $i_1$.

From the description of $L_n$ below Figure \ref{linkspic}, it is clear that $M_n$ may be described as

\[ B^3 - S \cup_{j_0} (S^2 \times I) - T \cup_{j_1} \hdots \cup_{j_1} (S^2 \times I) - T \cup_{j_0^{-1}} B^3 - \overline{S}. \]

We summarize the above discussion in the following lemma.

\begin{lemma}

$M_n$ is homeomorphic to
\[ M_S \cup_{i_0} M_T \cup_{i_1} \hdots \cup_{i_1} M_T \cup_{i_0^{-1}} \overline{M_S}, \]
by a homeomorphism that restricts on each component tangle complement to $h_S$ or $h_T$ as appropriate.

\end{lemma}

\begin{proof}[Proof of Theorem \ref{Ln}.]  First we claim that $M_S \cup_{i_0} M_T \cup_{i_1} \hdots \cup_{i_1} M_T$ (with $n$ copies of $M_T$) is isometric to the convex core of

\[ \left\langle\, G^{c^{2n}},H^{c^{2(n-1)}}, \hdots , H \, \right\rangle \]
by an isometry which restricts to $i_T$ on the rightmost copy of $M_T$.  We prove this by induction on $n$.

In what follows we denote by $\phi \co C(H) \rightarrow C(H^{c^{-2}})$ the isometry between convex cores induced by conjugation by $c^{-2}$.  Observe that $H^{c^{-2}} = \overline{H}$.  Since $H$ and $G$ both contain the subgroup $\Lambda \subset \mathrm{PSL}_2(\mathbb{R})$ corresponding to a convex core boundary component, it follows that $G \cap H^{c^{-2}} = \Lambda$ and $H \cap H^{c^{-2}} = \Lambda$.

Now for the case $n=1$, notice that $M_T$ is isometric to the convex core of $H^{c^{-2}}$ by $\phi \circ i_T$, which extends the isometry between $M_{T_0}$ and the convex core of $H_0$.  From the commutative diagram above, the isometry $i_0 \co F \rightarrow F_0$, between the surface $F = \bound M_S$ and the surface $F_0 \subset M_{T_0} \subset M_T$, is induced by identifying each with the quotient of $\mathcal{H}$ by $\Lambda = G \cap H_0 = G \cap H^{c^{-2}}$.  Hence by Lemma \ref{Maskit}, $M_S \cup_{i_0} M_{T}$ is isometric to the convex core of $\langle\, G, H^{c^{-2}} \,\rangle$.         Thus $M_S \cup_{i_0} M_{T}$ is also isometric to the convex core of $\langle \, G^{c^2},H \, \rangle$.  The restriction of this isometry to $M_T$ is $i_T$.

 For the induction step, Lemma \ref{Maskit} and the fact that $H \cap H^{c^{-2}} = \Lambda$ imply that $M_T \cup_{i_1} M_T$ is isometric to the convex core of $\langle\, H, H^{c^{-2}}\,\rangle$, by an isometry which restricts to $i_T$ and $\phi \circ i_T$ on the respective copies of $M_T$.  But by the induction hypothesis $i_T$ extends to an isometry of $M_S \cup_{i_0} M_T \cup_{i_1} \hdots \cup_{i_1} M_T$, with $n$ copies of $M_T$, with the convex core of $\langle\, G^{c^{2n}}, H^{c^{2(n-1)}}, \hdots, H \,\rangle$.  Thus there is an isometry from $M_S \cup_{i_0} M_T \cup_{i_1} \hdots \cup_{i_1} M_T$, with $n+1$ copies of $M_T$, to the convex core of $\langle\,  G^{c^{2n}}, H^{c^{2(n-1)}},\hdots,H,H^{c^{-2}} \, \rangle$ which restricts appropriately on all copies of $M_S$ and $M_T$.  By conjugating this group by $c^2$, we verify the claim.

The theorem now follows from the claim by applying Lemma \ref{Maskit}, using the same philosophy, to $\Gamma_n$, noting that $\overline{M_S}$ is (orientation--preserving) isometric to the convex core of $\overline{G}$ by the composition of $i_S$ with the isometry induced by reflection through  $\mathcal{H}$.  \end{proof}

\begin{cor1}  For any $n$, the trace field of $M_n$ is $\mathbb{Q}(i,\sqrt{2})$.  \end{cor1}

\begin{proof}  Our description of $\Gamma_n$ makes it clear that the trace field of $M_n$ is contained in the field generated by the traces of the group $\langle\, G, H, c\, \rangle$.  Since $\langle\, G, H, c\, \rangle$ is a subgroup of  $\mathrm{PSL}_2(\mathbb{Q}(i,\sqrt{2}))$ the trace field of $M_n$ is a subfield of $\mathbb{Q}(i,\sqrt{2})$.  On the other hand,

$$\tr(s)=1+i \AND \tr(h)=-i\sqrt{2}.$$ 
Since $\Gamma_n$ contains conjugates of both $s$ and $h$ we have that the trace field of $M_n$ is exactly $\Q(i, \sqrt{2})$.

\end{proof}

We use the {\it complex modulus} to show that the manifolds $M_n$ are pairwise not commensurable.

\begin{dfn}  Let $C$ be a rank two cusp of a hyperbolic manifold $M$ and take $\alpha, \beta \in \pi_1(M)$ to be peripheral elements corresponding to $C$ which generate a group of rank two.  Then $\alpha$ and $\beta$ can be simultaneously conjugated in $\mathrm{PSL}_2(\C)$ so that

$$\alpha \ =\ \begin{pmatrix} 1 & 1 \\ 0 & 1 \end{pmatrix}    \AND \beta \ =\ \begin{pmatrix} 1 & z \\ 0 & 1 \end{pmatrix}.$$
The complex modulus for $C$ is the equivalence class $m(C) = [z] \in \mathbb{C}/\mathrm{PGL}_2(\mathbb{Q}),$  where $\mathrm{PGL}_2(\mathbb{Q})$ acts on $\mathbb{C}$ by M\"obius transformations. 

\end{dfn}

\begin{lemma} Both cusps of $M_n$ have modulus $\left[ i(1+2n\sqrt{2}) \right].$\end{lemma}

\begin{proof}

Define 

\begin{align*}
& \gamma_1 \ = \ \inv{g}h\inv{(hg)} \\
& \gamma_2 \ = \ \mathrm{id} \\
& \gamma_3 \ = \ f\inv{h}fg \\
& \gamma_4 \ = \ \inv{g}fg
\end{align*}
and 
\[ \delta_j \ = \ (\gamma_j^{-1})^{c^2}\overline{\gamma_j}. \]  
Then $\delta_j \in H$ and $p_j^{\delta_j} \ = \ p_j^{c^2}$ for each $j$.  Let
\begin{align*}
& \lambda_1 \ = \ \overline{\inv{(srs)}} \, \delta_3^{-1} \, (\delta_3^{-1} )^{c^2} \ldots (\delta_3^{-1})^{c^{2(n-1)}} \, (srs)^{c^{2n}} \, \delta_1^{c^{2(n-1)}} \ldots \delta_1 \\
&  \lambda_2 \ = \ \overline{\inv{(sr\inv{s})}} \, \delta_4^{-1} \, (\delta_4^{-1} )^{c^2} \ldots (\delta_4^{-1})^{c^{2(n-1)}} \, (sr\inv{s})^{c^{2n}}. 
\end{align*}  
Then $\lambda_1, \, \lambda_2 \in \Gamma_n$ and $\lambda_1$ and $\lambda_2$ commute with $p_1$ and $p_2$, respectively.

Set $\tau =  \left( \begin{smallmatrix} 0 & i \\ i & 0 \end{smallmatrix} \right)$.  Then $p_1^\tau  = \left( \begin{smallmatrix}  1 & 1 \\ 0 & 1 \end{smallmatrix} \right)$ and, by a straightforward induction argument,
\[ \lambda_1^\tau \ = \  \begin{pmatrix} 1 & 2i(1+2n\sqrt{2}) \\ 0 & 1 \end{pmatrix}. \]
Therefore, $\langle p_1, \lambda_1 \rangle$ is a rank two subgroup for the knotted cusp $C'_n$ of $M_n$ and $m(C'_n) \ = \ \left[ i(1+2n\sqrt{2}) \right]$.

For the other cusp $C''_n$ of $M_n$, we set $\omega = \left( \begin{smallmatrix} i \sqrt{5}/5 & 0 \\ 0 & -i\sqrt{5} \end{smallmatrix} \right)$.  Then $p_2^\omega  = \left( \begin{smallmatrix} 1 & 1 \\ 0 & 1 \end{smallmatrix} \right)$ and
\[ \lambda_2^\omega \ = \ \begin{pmatrix} 1 & \frac{2i}{5}(1+2n\sqrt{2}) \\ 0 & 1\end{pmatrix}. \]
Therefore $m(C''_n) \ = \ \left[ i(1+2n\sqrt{2}) \right]$.

\end{proof}

The complex modulus of the cusp of a finite volume hyperbolic manifold is an invariant of the commensurability class of the Euclidean structure inherited by the cusp torus.  Since commensurable hyperbolic manifolds have commensurable cusps, the collection of moduli of the cusps of a finite volume hyperbolic manifold is a commensurability invariant.  This was first used by Thurston (\cite{Th}, chapter 6) to show that certain hyperbolic links are not commensurable.  In Thurston's examples, representatives for the cusp moduli of the incommensurable manifolds generate distinct number fields.  With our examples, representatives for the the complex modulus of the cusps all generate the trace field $\mathbb{Q}(i, \sqrt{2})$ so,  in order to prove the lemma below, we must use finer properties of of the $\mathrm{PSL}_2(\mathbb{Q})$--action on the trace field.

\begin{lem1}

$M_m$ and $M_n$ are not commensurable whenever $m \neq n.$

\end{lem1}

\begin{proof}  Suppose $M_m$ and $M_n$ are commensurable; then we have
\[ i(1+2n\sqrt{2}) = \frac{ai(1+2m\sqrt{2})+b}{ci(1+2m\sqrt{2})+d}, \]
where $ad-bc \neq 0$.  By clearing denominators, we may take $a,b,c,d \in \mathbb{Z}$.  The left-hand side above is imaginary, so setting the real part of the right-hand side above equal to zero yields

\[ bd - ac(1+8m^2) - ac(4m\sqrt{2}) = 0. \]

Since $\{1,\sqrt{2}\}$ is a basis of $\mathbb{Q}(i,\sqrt{2})$ as a vector space over $\mathbb{Q}$, we must have $ac = 0$ and therefore $bd = 0$ as well.  Since $ad - bc \neq 0$, either $b=c=0$ or $a=d=0$.  In the first case we obtain $i(1+2n\sqrt{2}) = \frac{a}{d}i(1+2m\sqrt{2})$.  This can only happen if $a= d$ and $m = n$.  In the second case, we have
\[ i(1+2n\sqrt{2}) = \frac{b}{ci(1+2m\sqrt{2})} = i\frac{b}{c(8n^2-1)}(1-2m\sqrt{2}). \]
This cannot occur, since the coefficients of 1 and $\sqrt{2}$ have the same sign on the left--hand side of the equation but opposite signs on the right.

\end{proof}

We have shown that each of the 2-component link complements $M_n$ represents a distinct commensurability class and has trace field $\mathbb{Q}(i,\sqrt{2})$.  To complete the proof of Theorem \ref{links}, it remains to produce links with $k+1$ components, where $k \geq 2$.

For each $n$, the two components of $L_n$ have linking number 0 and one of the two components is unknotted.  Hence, the preimage $L_n^{(k)}$ of $L_n$ under the $k$-fold branched cover of $S^3$ branched over the unknotted component of $L_n$ is a $(k+1)$-component link in $S^3$.  Its complement $M_n^{(k)}$ is a cyclic $k$-fold cover of $M_n$.  Therefore, for fixed $k$, the commensurability classes of the $M_n^{(k)}$ are exactly the commensurability classes of the $M_n$.  These examples complete the proof of Theorem \ref{links}.


\section{Nonintegral traces}  \label{sec:nonint}

\begin{dfn}  Let $\Gamma < \mathrm{SL}_2(\overline{\mathbb{Q}})$, where $\overline{\mathbb{Q}}$ is the algebraic closure of $\mathbb{Q}$ in $\mathbb{C}$.  $\Gamma$ has \textit{integral traces} if for all $\gamma \in \Gamma$, the trace of $\gamma$ is an algebraic integer.  Otherwise, $\Gamma$ has \textit{nonintegral traces}.  \end{dfn}

It is easy to see that, for $n \in \mathbb{N},$  $\Gamma_n$ has integral traces; indeed, the matrix entries of the elements of $\Gamma_n$ are all algebraic integers.  The incommensurability of the $M_n$ thus stands in stark contrast to the situation for finite--volume manifolds of the form $M=\mathbb{H}^3/\Gamma$, where $\Gamma < \mathrm{PSL}_2\left(\mathbb{Q}(\sqrt{-d})\right)$ for some positive square-free integer $d$.  \textit{All} such manifolds with the property that $\Gamma$ has integral traces are commensurable with $\mathbb{H}^3/\mathrm{PSL}_2(O_d)$, where $O_d$ is the ring of integers of $\mathbb{Q}(\sqrt{-d})$ (cf. \cite{MaR}, Theorem 8.3.2).  We give examples in Theorem \ref{L0} below.

In this section, we consider the family of links $L_n'$ produced by mutating $L_n$ along the 4-punctured sphere separating $S$ and $T$.  The group $\Gamma_n'$ which uniformizes the complement $M_n'$ of $L_n'$ contains an element with a nonintegral trace.  Since integrality of traces is a commensurability invariant (see eg. \cite{MaR}, Exercise 5.2, No. 1), this implies that $M_n$ and $M_n'$ are incommensurable.  Mutation along 4-punctured spheres is a much studied phenomenon, and hyperbolic manifolds related thereby can be difficult to distinguish.  Among other things, they have the same volume \cite{Ru} and invariant trace field \cite{NeR1}, and their A-polynomials have a common factor \cite{CL}.  In light of this, it may be surprising that in this case mutation produces a nonintegral trace.   However, there are other examples of this phenomenon in the literature, for instance the example of Figure 6.1 in \cite{LR}.  On the other hand, at the end of this section we give an example where mutation preserves the property of having integral traces.

It is well known that any hyperbolic 3-manifold homotopy equivalent to a 4-punctured sphere has a Klein 4-group of \textit{mutation} isometries, each of whose nontrivial elements acts on the set of punctures as an even permutation of order two (see eg. Theorem 2.2 of \cite{Ru}).  In what follows, we consider the totally geodesic boundary of $M_S$ and the mutation which exchanges the endpoints of one string of $S$ with the ends of the other.

Recall that the fundamental group $\Lambda$ of $\partial M_S$ is generated by parabolic elements $p_1,$ $p_2,$ and $p_3$ as defined in Section \ref{sec:S}, and the final conjugacy class of parabolic elements of $\Lambda$ is represented by $p_4=p_1p_2p_3^{-1}$.  Considering the words in $G$ representing the $p_i$'s, we find that $p_1$ and $p_3$ are conjugate in $G$, as are $p_2$ and $p_4.$  Hence $p_1$ and $p_3$ correspond to two ends of the same string of $S$ and $p_2$ and $p_4$ correspond to the two ends of the other string.

The fixed points of $p_1$, $p_2$, and $p_3$ are respectively $0$, $\infty$, and $5/3$.   Consider the element 
$$m \ = \  \begin{pmatrix} 0 & \sqrt{5} \\ -\frac{1}{\sqrt{5}} & 0 \end{pmatrix}$$
which exchanges $0$ and $\infty$.  It is easily verified that $m$ normalizes $\Lambda$ and acts on the generating parabolics as follows:

\begin{align*} 
 p_1^ m & = p_2 &
 &  p_2^m = p_1 &
 &  p_3^m =  p_4^{p_1^{-1}}. \end{align*}

\begin{figure}[ht]

\begin{center}

\includegraphics[height=1.25in]{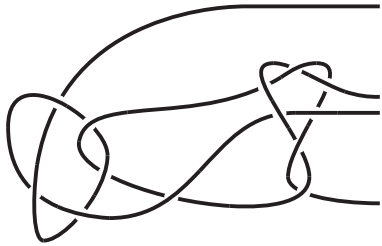}

\end{center}

\caption{The tangle formed by attaching $S$ and $T_0$ with the mutant isometry} \label{mutate}

\end{figure}

Thus the isometry induced by $m$ on the totally geodesic boundary of $M_S$ exchanges the ends of one string of $S$ with the ends of the other.  It follows from Lemma \ref{Maskit} that the manifold formed by gluing $M_S$ to $M_{T_0}$ by this isometry is isometric to $C(\langle\, G^m,H_0\, \rangle)$.  This manifold is homeomorphic to the complement in $B^3$ of the tangle pictured in Figure \ref{mutate}.  The group $\langle\, G^m,H_0\, \rangle$ contains the element
\[ msm^{-1}g \ = \ \begin{pmatrix} -1+\sqrt{2}+11i+i\sqrt{2} & 1-7\sqrt{2}-16i-2i\sqrt{2} \\ \frac{1}{5}(2-\sqrt{2}-21i-2i\sqrt{2}) & \frac{1}{5}(-2+12\sqrt{2}+31i+4i\sqrt{2}) \end{pmatrix} \]
whose trace has minimal polynomial 
$$ 25x^4-60x^3-106x^2+836x+13673.$$
Therefore $\mathrm{Tr}(msm^{-1}g)$ is not an algebraic integer.

\bigskip

Let $L_n'$ be the link obtained by performing this mutation along the totally geodesic 4-punctured sphere separating $S$ and $T$ in $M_n$.  Following the discussion of the description of the groups $\Gamma_n$, we conclude that the complement $M_n'$ of $L_n'$ is uniformized by

\[ \Gamma_n' = \left \langle\, (G^m)^{c^{2n}},H^{c^{2(n-1)}},\dots,H,\bar{G}\, \right \rangle. \]

Because traces are invariant under conjugation, $\Gamma_n'$ has a nonintegral trace for each $n \in \mathbb{N}$.  We summarize this in the theorem below.

\begin{thm}  For each $n$, $S^3 - L_n'$ is homeomorphic to $\mathbb{H}^3/\Gamma_n'$, and $\Gamma_n'$ has a nonintegral trace.  \end{thm}

Now consider the link we call $L_0$, formed by doubling the tangle $S$ without inserting any copies of $T$.  The complement of $L_0$ is uniformized by the group $\Gamma_0 = \langle\, G,\bar{G}\, \rangle$.  The complement of the link $L_0'$ resulting from the same mutation as above is uniformized by $\Gamma_0' = \langle\, G^m, \bar{G}\, \rangle$, which has generating set $\{r,msm^{-1},\bar{s}\}$.  The trace field of this group is generated over $\mathbb{Q}$ by the traces of the following elements (see \cite{MaR}, p. 125):
\[  r, \ msm^{-1}, \ \bar{s}, \ rmsm^{-1}, \ r\bar{s}, \ msm^{-1}\bar{s}, \ \mathrm{and} \ rmsm^{-1}\bar{s}.  \]
Of these words, the only ones which could possibly have nonintegral traces are the last two, but a computation reveals 
$$ \tr(msm^{-1}\bar{s}) \ = \ -10  \AND \tr(rmsm^{-1}\bar{s}) \ = \ -18+6i. $$
Thus, $\Gamma_0$ has integral traces.

\begin{thm}  The link $L_0$ and its mutant $L_0'$ each have hyperbolic complement commensurable with $\mathbb{H}^3/\mathrm{PSL}_2(O_1)$, where $O_1$ is the ring of integers of the number field $\mathbb{Q}(i)$.  \label{L0} \end{thm}

\begin{proof}  The description above shows that $\langle\,G,\overline{G}\,\rangle$ and $\langle\,G^m,\overline{G}\,\rangle$, which uniformize $L_0$ and $L_0'$, respectively, have traces in $O_1$.  Thus $L_0$ and $L_0'$ each satisfy the criterion mentioned at the beginning of this section, and are commensurable with $\mathbb{H}^3/\mathrm{PSL}_2(O_1)$.

\end{proof}


\section{One cusped examples} \label{sec:onecusp}

In this section, we use our techniques for constructing incommensurable manifolds with the same trace field to give examples with only one cusp.   As before, we produce manifolds $N$ and $N'$ as identification spaces of the ideal octahedron and cuboctahedron, respectively.  Again, these manifolds are constructed to have isometric totally geodesic boundaries.  This time we ensure that $\partial N$ has isometry group of full order, which allows copies of $N$ and $N'$ to be glued so that the resulting manifold has only one cusp.  The cost, as we will see, is that the $\mathbb{Z}/2\mathbb{Z}$ homology of this manifold is much larger than that of a knot complement.

We construct $N$ by identifying faces of the octahedron using the isometries
\[ a =  \begin{pmatrix} i & -i \\ 1-i & -1 \end{pmatrix} \AND  b = \begin{pmatrix} 2i & 2-i \\ i & 1-i \end{pmatrix}.  \]
With the octahedron positioned as in our earlier construction, $a$ takes $X_1$ to $X_2$ so that the ideal vertex $\frac{1+i}{2}$ maps to $\infty$, and $b$ takes $X_3$ to $X_4$ so that $1+i$ maps to $\infty$.

The identifications on the cuboctahedron are given by 
$$\begin{array}{rcl}

 x & = & \begin{pmatrix} i & -i \\ -i - \sqrt{2} & \sqrt{2} \end{pmatrix}  \\

 y & =  & \begin{pmatrix} 1-2i\sqrt{2} & 2+i\sqrt{2} \\ 1-i\sqrt{2} & 1+i\sqrt{2} \end{pmatrix} \\

 z & = & \begin{pmatrix} 4i & -3i-2\sqrt{2} \\ i-2\sqrt{2} & \sqrt{2}-3i \end{pmatrix}.  \end{array}$$
With the cuboctahedron as previously, $x$ takes $Y_2$ to $Y_1$ so that 0 maps to $-i\frac{\sqrt{2}}{2}$, $y$ takes $Y_3'$ to $Y_3$ so that $\frac{1}{3}(1-2i\sqrt{2})$ maps to $\infty$, and $z$ takes $Y_1'$ to $Y_2'$ so that $\frac{1}{2}(1-i\sqrt{2})$ maps to $\frac{1}{3}(2-2i\sqrt{2})$.  We refer to the resulting identification space as $N'.$

Take $G = \langle \, a,b \, \rangle$ and $H = \langle \, x,y,z \, \rangle$.  Then $N$ and $N'$ are isometric to the convex cores $C(G)$ and $C(H)$, respectively.  The totally geodesic boundary of $N$ is isometric to a totally geodesic boundary component of $N'$, which is reflected by the fact that $G$ and $H$ share a subgroup $\Lambda$ stabilizing the hyperplane $\mathcal{H}$ over $\mathbb{R}$.  $\Lambda$ is generated by the parabolic elements 

\[ \begin{array}{rcccccl}

 p_1 & = & aba & = & xyx & = &  \begin{pmatrix} 1&0 \\ -3&1 \end{pmatrix}  \\

 p_2 & = & b^2a^{-1} & = & yzx^{-1} & = & \begin{pmatrix} -1&-3 \\ 0&-1 \end{pmatrix}  \\

 p_3 & = & a^{-2}b^{-1} & = & x^{-2}y^{-1} & = & \begin{pmatrix} -2&3 \\ -3&4 \end{pmatrix}. \end{array} \]
We have $p_1(0)=0,$ $p_2(\infty)=\infty,$ and $p_3(1)=1.$  The element 
$$p_4 \ = \ p_3 p_2 p_1 \ = \ a^{-2}b^2a \ = \ x^{-2}zyx \ \in \Lambda$$ 
fixes $1/2$ and represents the fourth cusp of the boundary.  Note that $p_1$ and $p_3$ are conjugate in $G$, as are $p_2$ and $p_4.$  Also, $p_1$ and $p_3$ are conjugate in $H$.  Since $p_1$ and $p_3$ are conjugate in both groups, the convex core of $\mathbb{H}^3/\langle  G, H  \rangle$ will have a closed cusp.  Thus in order to build manifolds with a single cusp, we identify boundary components of $N$ and $N'$ by a nontrivial isometry.

\begin{remark}

Notice that $\Lambda$ is the kernel of the map $\mathrm{PSL}_2(\mathbb{Z}) \rightarrow \mathrm{PSL}_2(\mathbb{Z}/3\mathbb{Z})$ given by entrywise reduction.  Thus, $\mathrm{PSL}_2(\mathbb{Z})$ acts on $\Lambda$ by conjugation inducing isometries of the 4-punctured sphere uniformized by $\Lambda$.  This contrasts the behavior of the totally geodesic 4-punctured spheres in the links $L_n$, whose groups are not normal in $\mathrm{PSL}_2(\mathbb{Z})$ and have normalizers which include elements not in $\mathrm{PSL}_2(\mathbb{Z})$.  Therefore every manifold obtained by gluing $N$ to $N'$ by an isometry of the boundary has integral traces.

\end{remark}

The parabolic element $m\,=\,\left( \begin{smallmatrix} 1&1 \\ 0&1 \end{smallmatrix} \right)\, \in\, \mathrm{PSL}_2(\mathbb{Z})$ conjugates $\Lambda$ to itself with action given by
 \begin{align*} & p_1^m = p_3 & & p_2^m = p_2 & & p_3^m = (p_4^{-1})^{p_3^{-1}}. \end{align*}
This induces an isometry of $\mathcal{H}/\Lambda$ which acts fixing the cusp corresponding to $p_2$ and cyclically permuting the other three.

The other boundary component of $N'$ is a quotient of the geodesic hyperplane over $\mathbb{R}-i\sqrt{2}$.  The generators for the stabilizer in $H$ of this hyperplane are $q_1 = xy^{-1}z^{-1}$, $q_2 = yzx^{-1},$ and $q_3 = z^2y^{-1}.$  Then 
$$q_1\left(-i\sqrt{2}\right)\ =\ -i\sqrt{2}, \quad q_2\left(\infty\right)\ =\ \infty, \quad \mathrm{and} \ \  q_3\left(1-i\sqrt{2}\right)\ =\ 1-i\sqrt{2}.$$  
Each $q_i$ corresponds to $p_i$ in the following way: take
\[ c \ = \ \begin{pmatrix} 1 & i\sqrt{2} \\ 0 & 1 \end{pmatrix}. \]
Then for each $i$, $cq_ic^{-1} = p_i$.  (Set $q_4 = q_3 q_2 q_1 = z^3y^{-1}z^{-1}$).  In particular, we note that the two totally geodesic boundary components of $N'$ are isometric.

\begin{thm} For each $n \in \mathbb{N}$, $N_n = \mathbb{H}^3/\Gamma_n$ has one cusp, where $\Gamma_n$ is defined by 

\[ \Gamma_n \ = \ \left\langle\ (G^m)^{c^n}, H^{c^n}, \hdots, H^c,\overline{G}\ \right\rangle. \]  \end{thm}

\begin{proof}  By a line of argument similar to the proof of Theorem \ref{Ln}, $\mathbb{H}^3/\Gamma_n$ is isometric to the hyperbolic manifold obtained by stacking $n$ copies of $N'$ end--to--end, capping off on one end by $N$ and on the other by $\overline{N}$.  The isometry of boundary components by which $N$ is glued to the first copy of $N'$ is induced by $m$ as above.  The isometries by which copies of $N'$ are glued to each other take the cusp corresponding to each $q_i$ to the cusp corresponding to $p_i$, as does the isometry taking the boundary of $\overline{N}$ to the ``back'' boundary of the final copy of $N'$.

\vspace{.75in}

\begin{figure}[ht]

\begin{picture}(100,100)

\put(10,15) {\includegraphics[width=4.5in]{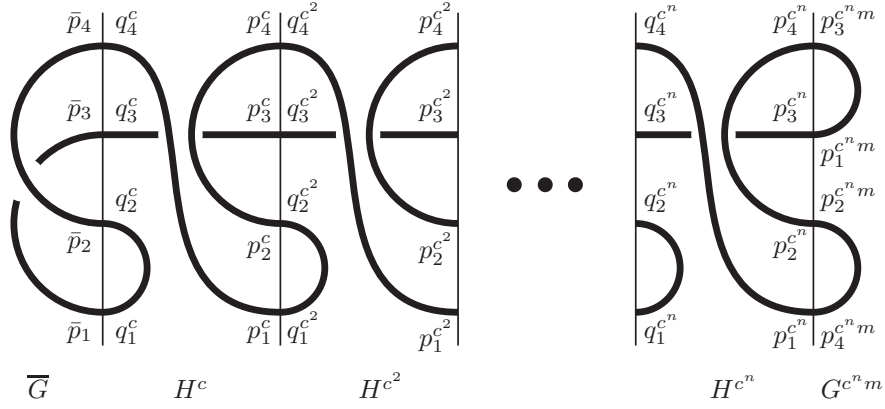}}

\put(35,18){$\bar{p}_1$}  \put(53,18){$q_1^c$}

\put(35, 52){$\bar{p}_2$} \put(53, 68){$q_2^c$}

\put(35, 102){$\bar{p}_3$} \put(53, 102){$q_3^c$}

\put(35, 135){$\bar{p}_4$} \put(53,135){$q_4^c$}

\put(103,18){$p_1^c$}  \put(118,18){$q_1^{c^2}$}

\put(103, 50){$p_2^c$} \put(118, 68){$q_2^{c^2}$}

\put(103, 102){$p_3^c$} \put(118, 102){$q_3^{c^2}$}

\put(103, 135){$p_4^c$} \put(118,135){$q_4^{c^2}$}

\put(168,15){$p_1^{c^2}$}  

\put(168, 48){$p_2^{c^2}$} 

\put(168, 102){$p_3^{c^2}$} 

\put(168, 135){$p_4^{c^2}$} 

\put(253,18){$q_1^{c^n}$}

\put(253, 68){$q_2^{c^n}$}

\put(253, 102){$q_3^{c^n}$}

\put(253,135){$q_4^{c^n}$}

\put(302,16){$p_1^{c^n}$}  

\put(302, 50){$p_2^{c^n}$} 

\put(302, 102){$p_3^{c^n}$} 

\put(302, 135){$p_4^{c^n}$} 

\put(320,16){$p_4^{c^n m}$}

\put(320,68){$p_2^{c^n m}$}

\put(320,85){$p_1^{c^n m}$}

\put(320, 135){$p_3^{c^n m}$} 

\put(20,-5){$\overline{G}$}

\put(75,-5){$H^c$}

\put(145,-5){$H^{c^2}$}

\put(278,-5){$H^{c^n}$}

\put(320,-5){$G^{c^n m}$}

\end{picture}

\caption{A schematic for $N_n$}  \label{schematic}

\end{figure}


In $H$, $q_1$ is conjugate to $p_4$, whereas $p_2$ and $q_2$ are conjugate and so are $q_3$ and $q_4$.  The schematic of Figure \ref{schematic} incorporates this information as well as the conjugacy information from $G$ to show how the annular cusps of $N$ and $N'$ are identified along their boundaries in $N_n$.  From this it is evident that $\mathbb{H}^3/\Gamma_n$ has only one cusp.

\end{proof}

As before, we can use the complex modulus to show that the manifolds $N_n$ are pairwise not commensurable.

\begin{lemma} The complex modulus for the cusp of $N_n$ is $\left[ 1+4i(n\sqrt{2}+1) \right].$  \end{lemma}

\begin{proof} Let 

\begin{align*}
\alpha \ &= \ y \inv{z} \inv{c} \inv{x} c \inv{y} \\
\beta \ &= \ \inv{a} \inv{m} \inv{x} ma \inv{b} \\
\gamma \ &= \ \overline{b\inv{a}} \, c (x\inv{y}) \inv{c} \, \overline{ba}.
\end{align*}
If we set $\mu  =   \overline{aba}$ and
\[ \lambda_n \ = \ (\inv{y})^c \,  \alpha^{c^2} \cdots \alpha^{c^n} \,  \beta^{c^nm}  \, \gamma \]
then $\mu$ and $\lambda_n$ are elements of the group $\Gamma_n$.


Let $\tau =  \left( \begin{smallmatrix} 0 & \sqrt{3}/3 \\ -\sqrt{3} & 0 \end{smallmatrix} \right)$.  Then  $\mu^\tau  = \left( \begin{smallmatrix} 1 & 1 \\ 0 & 1 \end{smallmatrix}\right)$ and
 \[ \tau \lambda_n \inv{\tau} \ = \ \begin{pmatrix} 1 & \frac{1}{3}\left( 1+4i(n\sqrt{2}+1) \right) \\ 0 & 1\end{pmatrix}. \]
Therefore, the complex modulus for the manifold $N_n$ is $\left[ 1+4i(n\sqrt{2}+1) \right].$

\end{proof}

\begin{lemma}

$N_m$ and $N_n$ are not commensurable whenever $m \neq n.$

\end{lemma}

\begin{proof}

As in the two component link case, it is easy to check that if $\phi$ is a M\"obius transformation in $\mathrm{PGL}_2(\mathbb{Q})$ with $\phi \left( 1+4i(m\sqrt{2}+1)  \right) =  1+4i(n\sqrt{2}+1) .$ then $m=n.$

\end{proof}

Thus, our construction gives infinitely many commensurability classes of one-cusped hyperbolic manifolds, and it is easily seen that all of the $N_n$ have trace field $\mathbb{Q}(i,\sqrt{2})$.  For hyperbolic manifolds which are not link complements (and we will see below that these are not), the trace field may fail to be a commensurability invariant.  For general hyperbolic manifolds $M = \mathbb{H}^3/\Gamma$, the \textit{invariant trace field} is a commensurability invariant.  This is defined to be the trace field of the finite--index subgroup
\[ \Gamma^{(2)} = \{\ \gamma^2\ |\ \gamma \in \Gamma\ \}. \]
See \cite{MaR}, \S 3.3 for a thorough discussion of these invariants.  In the current situation, the trace field and invariant trace field coincide.

\begin{lemma}  For any $n \in \mathbb{N}$, $N_n$ has invariant trace field $\mathbb{Q}(i,\sqrt{2})$.  \end{lemma}

\begin{proof}  Inspection of the traces of $a$ and $b$ shows that the trace field of $N$ is $\mathbb{Q}(i)$.  The invariant trace field is a subfield which is not properly contained in $\mathbb{R}$, since $\langle\,a,b\,\rangle$ is not a Fuchsian group (cf. \cite{MaR}, Theorem 3.3.7).  Thus the invariant trace field of $N$ is also $\mathbb{Q}(i)$.  Similarly, one finds that the invariant trace field of $N'$ is $\mathbb{Q}(i\sqrt{2})$.  The proof is completed using induction and a theorem of Neumann-Reid \cite{NeR1}, which asserts that the invariant trace field of a Kleinian group obtained as a free product with amalgamation $A *_C B$ is the compositum of the invariant trace fields of $A$ and $B$.  \end{proof}

This completes the proof of Theorem \ref{onecusped} stated in the introduction.  Given that the motivating question concerns knot complements, we would like to determine whether or not the $N_n$ are knot complements in $S^3$.  The following proposition shows that they are, in fact, quite far from this.

\begin{prop1} 

$$\mathrm{H}_1(N_n;\mathbb{Z}/2\mathbb{Z}) \ \cong \ (\mathbb{Z}/2\mathbb{Z})^{n+1}.$$ 

\end{prop1}

\begin{proof} For $i \in \{1,\hdots,n\}$, define $x_i = ,x^{c^i}$, $y_i = y^{c^i}$, and $z_i = z^{c^i}$.  Let $a_n = a^{c^nm}$, $b_n = b^{c^nm}$, $a_0 = \bar{a}$, and $b_0 = \bar{b}$.  Our description of $\Gamma_n$ shows that it has a generating set 

\[ \left(\bigcup_{i = 1}^n \{ x_i,y_i,z_i \}\right) \cup \{a_n,b_n,\bar{a},\bar{b}\}. \]

Repeated applications of the Klein-Maskit combination theorem give a presentation for $\Gamma_n$ with ``middle relations'' \begin{align*}
  x_iy_ix_i & =  x_{i+1}y_{i+1}^{-1}z_{i+1}^{-1} \\
  y_iz_ix_i^{-1} & =  y_{i+1}z_{i+1}x_{i+1}^{-1} \\
  x_i^{-2}y_i^{-1} & =  z_{i+1}^2y_{i+1}^{-1}  \end{align*}
for $1 \leq i < n.$  We also get ``cap relations'' \begin{align*}
  a_nb_na_n & = x_n^{-2}y_n^{-1} &  a_0 b_0 a_0 & = x_1y_1^{-1}z_1^{-1} \\
  b_n^2a_n^{-1} & = y_nz_nx_n^{-1}  & b_0^2a_0^{-1} & = y_1z_1x_1^{-1} \\ 
  a_n^{-2}b_n^{-1} & = y_nx_ny_n^{-1}z_n^{-1}y_n^{-1}  & a_0^{-2} b_0^{-1} &= z_1^2 y_1^{-1}. \end{align*}
Abelianize these relations and reduce modulo 2 to get a presentation for the $\mathbb{Z}/2\mathbb{Z}$ homology of $N_n$.  The middle relations become \begin{align*} & y_i = x_{i+1}y_{i+1}z_{i+1} & & x_i y_i z_i = x_{i+1}y_{i+1}z_{i+1} & & y_i = y_{i+1} \end{align*}
which reduce to $x_i = z_i$, $x_{i+1}=z_{i+1}$, and $y_i = y_{i+1}$.  The cap relations become \begin{align*}
 & b_n = y_n & & b_0 = x_1 y_1 z_1 \\
 & a_n = x_n y_n z_n & & a_0 = x_1 y_1 z_1 \\
 & b_n = x_n y_n z_n & & b_0 = y_1.  \end{align*}
These reduce to $a_n = b_n = y_n$, $x_n = z_n$ and $a_0 = b_0 = y_1$, $x_1 = z_1$.  A basis for the $(\mathbb{Z}/2\mathbb{Z})$-vector space $\mathrm{H}_1(N_n;\mathbb{Z}/2\mathbb{Z})$ is therefore given by $\{x_1,\hdots,x_n,a_0\}$.

\end{proof}

\section{Further questions}

The motivating question remains unanswered.  Namely, do there exist infinitely many knot complements with trace field of bounded degree?  There does not seem to be an \textit{a priori} reason that our methods cannot be applied to construct knot complements.  However, at this point we do not have a starting place for our construction.  We have the following question (cf. the last sentence of \cite{Adams}).

\begin{question} Does there exist a hyperbolic knot complement in $S^3$ containing an embedded separating totally geodesic surface?  \end{question}

Related questions have been addressed in the literature.  Menasco-Reid conjectured that no hyperbolic knot complement contains a \textit{closed} embedded totally geodesic surface \cite{MeR}, and proved this conjecture for alternating knots.  It is known to be true for various other classes of knots, but false for links.  Indeed, Leininger constructed a hyperbolic $n$ component link containing a closed embedded totally geodesic surface for each $n \geq 2$ \cite{Le}.  Recently, Adams et. al. have constructed several infinite families of hyperbolic knot complements with totally geodesic Seifert surfaces \cite{A1}, \cite{A2}.

In another direction, there are many open questions about trace fields of hyperbolic 3-manifolds.  It is not known in general which number fields occur as trace fields; indeed, no number field $k \not\subset \mathbb{R}$ is known \textit{not} to arise this way.  This seems an interesting problem.  Our work suggests the following further question.

\begin{question} If a number field $k$ is the invariant trace field of a hyperbolic manifold, do there exist infinitely many commensurability classes of hyperbolic manifolds with invariant trace field $k$? 

\end{question}

A number field with one complex place admits infinitely many commensurability classes of compact arithmetic manifolds, distinguished by their \textit{invariant quaternion algebras}.  (See \cite{MaR} for an overview of arithmetic 3-manifolds and a definition of the invariant quaternion algebra.)   Note that the field $\mathbb{Q}(i, \sqrt{2})$ can be described as the compositum of  two fields $k_1$ and $k_2$ where each $k_j$ has one complex place and $k_1 \cap \mathbb{R} = k_2 \cap \mathbb{R}$.  Alan Reid has shown us that each field with this property admits infinitely many commensurability classes of {\it compact} 3-manifolds as well.  These may be constructed by gluing pieces of arithmetic manifolds along a totally geodesic surface so that their commensurability classes are distinguished by their invariant quaternion algebras.  Many fields arising as invariant trace fields of closed manifolds, for instance in the Snap census (cf. \cite{CGHN}), do not satisfy the criterion above, and for these fields we do not know the answer to this question.

Finite volume noncompact hyperbolic 3-manifolds with the same trace field cannot be distinguished by their invariant quaternion algebras (cf. \cite{MaR}, Theorem 3.3.8).  However, Ian Agol has pointed out to us that his methods of constructing hyperbolic manifolds with short systoles \cite{Agol} may be applied in dimension three to produce infinitely many commensurability classes of \textit{noncompact} hyperbolic manifolds with invariant trace field $\mathbb{Q}(\sqrt{-d})$, for any square-free positive integer $d$.  We know of no other noncompact examples but those of this paper.

\appendix

\section*{Appendix: Proof of Lemma \ref{convexcore}}

Following Morgan \cite{Mo}, we define a \textit{pared manifold} to be a pair $(M,P)$, where $M$ is a compact, orientable, irreducible 3-manifold with nonempty boundary which is not a 3-ball, and $P \subseteq \partial M$ is the union of a collection of disjoint incompressible annuli and tori satisfying the following properties: \begin{itemize}

  \item  Every noncyclic abelian subgroup of $\pi_1 M$ is conjugate into the fundamental group of a component of $P$.

  \item  Every map $\phi\co (S^1\times I, S^1\times \partial I) \rightarrow (M,P)$ which induces an injection on fundamental groups is homotopic as a map of pairs to a map $\psi$ such that $\psi(S^1\times I) \subset P$.

\end{itemize}

This definition is intended to capture the topology of the compact manifold obtained by truncating the cusps of the convex core of a geometrically finite hyperbolic 3-manifold by open horoball neighborhoods.  Indeed, Corollary 6.10 of \cite{Mo} asserts that if $(M,P)$ is obtained in this way, where $P$ consists of the collection of boundaries of the truncating horoball neighborhoods, then $(M,P)$ is a pared manifold.

Lemma \ref{convexcore} from the body of this paper asserts that if $(M,P)$ has the pared homotopy type of a geometrically finite hyperbolic manifold $\mathbb{H}^3/\Gamma$ where $\Gamma$ is not Fuchsian and $\partial C(\Gamma)$ is totally geodesic, then $M-P$ is homeomorphic to $C(\Gamma)$.  The key point of the proof is that the geometric conditions on $\Gamma$ ensure that $(M,P)$ is an acylindrical pared manifold.  Then Johannson's Theorem \cite{Jo}, that pared homotopy equivalences between acylindrical pared manifolds are homotopic to pared homeomorphisms, applies.  We expand on this below.

It is worth noting that Lemma \ref{convexcore} fails in more general circumstances.  Canary-McCullough give examples of this phenomenon in \cite{CM}, where for instance they describe homotopy equivalent non-Fuchsian geometrically finite manifolds with incompressible convex core boundary which are not homeomorphic (Example 1.4.5).  Their memoir \cite{CM} is devoted to understanding the ways in which homotopy equivalences of hyperbolic 3-manifolds can fail to be homeomorphic to homeomorphisms, and Lemma \ref{convexcore} follows quickly from results therein.

The treatment of Canary-McCullough itself uses the theory of {\it characteristic submanifolds} of manifolds with \textit{boundary pattern} developed in \cite{Jo}.  The characteristic submanifold of a manifold with boundary pattern is a maximal collection of disjoint codimension--zero submanifolds, each of which is an interval bundle or Seifert--fibered space embedded reasonably with respect to the boundary pattern.  Rather than attempting to establish all of the notation necessary to define this formally, we refer the interested reader to \cite{Jo} and \cite{CM}.  Here we simply transcribe the relevant theorem of \cite{CM}, which places strong restrictions on the topology of the characteristic submanifold of a pared manifold whose boundary pattern is determined by the pared locus.

For the purposes of Lemma \ref{convexcore} we exclude from consideration certain pared manifolds which never arise from convex cores of geometrically finite hyperbolic 3-manifolds.  We say $(M,P)$ is \textit{elementary} if it is homeomorphic to one of $(T^2\times I, T^2\times\{0\})$, $(A^2 \times I,A^2 \times \{0\})$, or $(A^2 \times I ,\emptyset)$, where $T^2$ and $A^2$ denote the torus and annulus, respectively;  otherwise $(M,P)$ is nonelementary.  Define $\partial_0 M := \overline{M-P}$.  We say an annulus properly embedded in $M-P$ is \textit{essential} in $(M,P)$ if it is incompressible and boundary--incompressible in $M-P$.  For a codimension--0 submanifold $V$ embedded in $M$, we denote by $\mathrm{Fr}(V)$ the \textit{frontier} of $V$ (that is, its topological boundary in $M$), and note that $\mathrm{Fr}(V) = \overline{\partial V - (V \cap \partial M)}$.  With notation thus established, the following theorem combines the definition of the characteristic submanifold with Theorem 5.3.4 of \cite{CM}.

\begin{thm1} Let $(M,P)$ be a nonelementary pared manifold with $\partial_0 M$ incompressible.   Select the fibering of the characteristic submanifold so that no component is an $I$--bundle over an annulus or M\"obius band.  \begin{enumerate}

  \item  Suppose $V$ is a component of the characteristic submanifold which is an $I$--bundle over a surface $B$.  Then each component of the associated $\partial I$--bundle is contained in $\partial_0 M$, each component of the associated $I$--bundle over $\partial B$ is either a component of $P$ or a properly embedded essential annulus, and $B$ has negative Euler characteristic.

  \item  Suppose $V$ is a Seifert fibered component of the characteristic submanifold.  Then $V$ is homeomorphic either to $T^2 \times I$ or to a solid torus.  If $V$ is homeomorphic to $T^2 \times I$, then one of its boundary components lies in $P$, the other components of $V \cap \partial M$ are annuli in $\partial_0 M$, and all components of $\mathrm{Fr}(V)$ are properly embedded essential annuli.  If $V$ is a solid torus, then $V \cap \partial M$ has at least one component, each an annulus either containing a component of $P$ or contained in $\partial_0 M$.  The components of $\mathrm{Fr}(V)$ are properly embedded essential annuli.  \end{enumerate} 

The characteristic submanifold contains regular neighborhoods of all components of $P$. \label{CM} \end{thm1}

The key claim in the proof of Lemma \ref{convexcore} is a further restriction on the characteristic submanifold of $(M,P)$, in the case that $M$ is obtained from the convex core of a non-Fuchsian geometrically finite manifold with totally geodesic convex core boundary by removing horoball neighborhoods of the cusps.  $P$ is the union of the boundaries of these neighborhoods.

\begin{claim}  $(M,P)$ as above is nonelementary, and $\partial_0 M$ is incompressible.  The characteristic submanifold of $(M, P)$ consists only of (Seifert fibered) regular neighborhoods of the components of $P$, each of whose boundary has a unique component of intersection with $\partial M$.  

\end{claim}

We prove the claim below, but assuming it for now, the proof of Lemma 3 proceeds as follows.  A representation as given in the statement of the lemma induces a pared homotopy equivalence between $(M,P)$ and the pared manifold $(N,Q)$ obtained by truncating $C(\Gamma)$ with open horoball neighborhoods.  Since $C(\Gamma)$ has totally geodesic convex core boundary, $(N,Q)$ is as described by the claim; hence $(M,P)$ is as well (see Theorem 2.11.1 of \cite{CM}, for example).  Johansson's Classification Theorem (cf. \cite{CM}, Theorem 2.9.10) implies that the original pared homotopy equivalence is homotopic to one which maps the complement of the characteristic submanifold of $(M,P)$ homeomorphically to the complement of the characteristic submanifold of $(N,Q)$.  It follows from the claim that these are homeomorphic to $M-P$ and $N-Q$, respectively, and the lemma follows.

\begin{proof}[Proof of claim]  As was mentioned above, the elementary pared manifolds do not arise from geometrically finite hyperbolic manifolds.  Since $(M, P)$ is obtained from the convex core of a geometrically finite manifold with totally geodesic convex core boundary, the following are known not to occur:

\begin{enumerate}

\item  A compressing disk for $\partial_0 M$.

(By definition $\partial_0 M$ lifts to a geodesic hyperplane in $\mathbb{H}^3$, hence the induced map $\pi_1 \partial M_0 \rightarrow \pi_1 M$ is injective.) 

\item  An \textit{accidental parabolic}; that is, an incompressible annulus properly embedded in $M$ with one boundary component in $P$ and one in $\partial_0 M$, which is not parallel to $P$.

(Every essential curve on $\partial_0 M$ that is not boundary-parallel is homotopic to a geodesic, but an element of $\pi_1(M)$ corresponding to an accidental parabolic has translation length $0$.)

\item  A \textit{cylinder}; that is, a properly embedded essential annulus in $M-P$.

(The double $DM$ of $M$ across $\partial_0 M$ is a hyperbolic manifold and the double of a cylinder in $M$ is an essential torus in $DM$.)

\end{enumerate}

We show that if the characteristic manifold has any components other than those listed in the claim then at least one of the above facts cannot hold.

For a component $V$ of the characteristic submanifold which is an $I$--bundle over a surface $B$, at least one component of the associated $I$--bundle over $\partial B$ must be properly embedded, since otherwise we would have $M= V$ and it is well known that an $I$--bundle over a surface does not admit a hyperbolic structure with totally geodesic convex core boundary unless the convex core is a Fuchsian surface.  But this annulus violates fact 2 or 3.  Thus there are no $I$--bundle components of the characteristic submanifold.

If $V$ is a Seifert fibered component of the characteristic submanifold homeomorphic to $T^2 \times I$, then one component of $\partial V$ is a torus $P_1 \subset P$, and all other components of $\partial V \cap \partial M$ are annulli in $\bound_0 M$.  If this second class is nonempty, then each component of $\mathrm{Fr}(V)$ is an essential annulus properly embedded in 

$M-P$.  This is not possible by fact 3, so $\partial V \cap \partial M$ consists only of $P_1$ and $V$ is a regular neighborhood of $P_1$.

If $V$ is a solid torus and $V \cap \partial M$ contains a component of $P$, then a similar argument shows that this is the unique component of $\partial V \cap \partial M$, so in this case $V$ is a regular neighborhood of an annular component of $P$.  If on the other hand $V \cap \partial M$ does \textit{not} contain any components of $P$, then it has at least two components, for otherwise a meridional disk of $V$ determines a boundary compression of the annulus $\mathrm{Fr}(V)$ in $M-P$.  But then any component of $\mathrm{Fr}(V)$ violates fact 3.  

\end{proof}


\bibliographystyle{plain}

\bibliography{commensurability}

\begin{thebibliography}{10}

\bibitem{Adams}
C.~Adams.
\newblock Hyperbolic knots.
\newblock In {\em Handbook of knot theory}, pages 1--18. Elsevier B. V.,
  Amsterdam, 2005.

\bibitem{A2}
C.~Adams, H.~Bennett, C.~Davis, M.~Jennings, J.~Novak, N.~Perry, and
  E.~Schoenfeld.
\newblock Totally geodesic {S}eifert surfaces in hyperbolic knot and link
  complements {II}.
\newblock preprint: arXiv.math.GT/0411358, November 2004.

\bibitem{A1}
C.~Adams and E.~Schoenfeld.
\newblock Totally geodesic {S}eifert surfaces in hyperbolic knot and link
  complements. {I}.
\newblock {\em Geom. Dedicata}, 116:237--247, 2005.

\bibitem{Agol}
I.~Agol.
\newblock Systoles of hyperbolic 4-manifolds, December 2006.

\bibitem{AR}
I.~R. Aitchison and J.~H. Rubinstein.
\newblock Combinatorial cubings, cusps, and the dodecahedral knots.
\newblock In {\em Topology '90 (Columbus, OH, 1990)}, volume~1 of {\em Ohio
  State Univ. Math. Res. Inst. Publ.}, pages 17--26. de Gruyter, Berlin, 1992.

\bibitem{Ba}
H.~Bass.
\newblock Finitely generated subgroups of {${\rm GL}\sb{2}$}.
\newblock In {\em The Smith conjecture (New York, 1979)}, volume 112 of {\em
  Pure Appl. Math.}, pages 127--136. Academic Press, Orlando, FL, 1984.

\bibitem{CM}
R.~D. Canary and D.~McCullough.
\newblock Homotopy equivalences of 3-manifolds and deformation theory of
  {K}leinian groups.
\newblock {\em Mem. Amer. Math. Soc.}, 172(812):xii+218, 2004.

\bibitem{CL}
D.~Cooper and D.~D. Long.
\newblock Remarks on the {$A$}-polynomial of a knot.
\newblock {\em J. Knot Theory Ramifications}, 5(5):609--628, 1996.

\bibitem{CGHN}
D.~Coulson, O.~A. Goodman, C.~D. Hodgson, and Walter~D. Neumann.
\newblock Computing arithmetic invariants of 3-manifolds.
\newblock {\em Experiment. Math.}, 9(1):127--152, 2000.

\bibitem{Jo}
K.~Johannson.
\newblock {\em Homotopy equivalences of {$3$}-manifolds with boundaries},
  volume 761 of {\em Lecture Notes in Mathematics}.
\newblock Springer, Berlin, 1979.

\bibitem{Ko1}
S.~Kojima.
\newblock Polyhedral decomposition of hyperbolic manifolds with boundary.
\newblock In {\em Proc. Workshops Pure. Math.}, volume~10, pages 37--57. 1990.

\bibitem{Ko2}
S.~Kojima.
\newblock Geometry of hyperbolic {$3$}-manifolds with boundary.
\newblock {\em Kodai Math. J.}, 17(3):530--537, 1994.
\newblock Workshop on Geometry and Topology (Hanoi, 1993).

\bibitem{Le}
C.~J. Leininger.
\newblock Small curvature surfaces in hyperbolic 3-manifolds.
\newblock {\em J. Knot Theory Ramifications}, 15(3):379--411, 2006.

\bibitem{LR}
D.~D. Long and A.~W. Reid.
\newblock Pseudomodular surfaces.
\newblock {\em J. Reine Angew. Math.}, 552:77--100, 2002.

\bibitem{MaR}
C.~Maclachlan and A.~W. Reid.
\newblock {\em The arithmetic of hyperbolic 3-manifolds}, volume 219 of {\em
  Graduate Texts in Mathematics}.
\newblock Springer-Verlag, New York, 2003.

\bibitem{Mask}
B.~Maskit.
\newblock On {K}lein's combination theorem. {III}.
\newblock In {\em Advances in the Theory of Riemann Surfaces (Proc. Conf.,
  Stony Brook, N.Y., 1969)}, Ann. of Math. Studies, No. 66, pages 297--316.
  Princeton Univ. Press, Princeton, N.J., 1971.

\bibitem{MeR}
W.~Menasco and A.~W. Reid.
\newblock Totally geodesic surfaces in hyperbolic link complements.
\newblock In {\em Topology '90 (Columbus, OH, 1990)}, volume~1 of {\em Ohio
  State Univ. Math. Res. Inst. Publ.}, pages 215--226. de Gruyter, Berlin,
  1992.

\bibitem{Mo}
J.~W. Morgan.
\newblock On {T}hurston's uniformization theorem for three-dimensional
  manifolds.
\newblock In {\em The Smith conjecture (New York, 1979)}, volume 112 of {\em
  Pure Appl. Math.}, pages 37--125. Academic Press, Orlando, FL, 1984.

\bibitem{NeR1}
W.~D. Neumann and A.~W. Reid.
\newblock Amalgamation and the invariant trace field of a {K}leinian group.
\newblock {\em Math. Proc. Cambridge Philos. Soc.}, 109(3):509--515, 1991.

\bibitem{NeR2}
W.~D. Neumann and A.~W. Reid.
\newblock Arithmetic of hyperbolic manifolds.
\newblock In {\em Topology '90 (Columbus, OH, 1990)}, volume~1 of {\em Ohio
  State Univ. Math. Res. Inst. Publ.}, pages 273--310. de Gruyter, Berlin,
  1992.

\bibitem{ZP}
L.~Paoluzzi and B.~Zimmermann.
\newblock On a class of hyperbolic {$3$}-manifolds and groups with one defining
  relation.
\newblock {\em Geom. Dedicata}, 60(2):113--123, 1996.

\bibitem{Ra}
J.~G. Ratcliffe.
\newblock {\em Foundations of hyperbolic manifolds}, volume 149 of {\em
  Graduate Texts in Mathematics}.
\newblock Springer-Verlag, New York, 1994.

\bibitem{Reid}
A.~W. Reid.
\newblock Arithmeticity of knot complements.
\newblock {\em J. London Math. Soc. (2)}, 43(1):171--184, 1991.

\bibitem{Ru}
D.~Ruberman.
\newblock Mutation and volumes of knots in {$S\sp 3$}.
\newblock {\em Invent. Math.}, 90(1):189--215, 1987.

\bibitem{Th}
W.~P. Thurston.
\newblock The geometry and topology of 3-manifolds.
\newblock mimeographed lecture notes, 1979.

\end{thebibliography}


\end{document}